\documentclass[12pt]{article}

\usepackage[english]{babel}
\usepackage[fixlanguage]{babelbib}
\usepackage[T1]{fontenc}
\usepackage[a4paper]{geometry}
\usepackage{amsmath}
\usepackage{amsthm}
\usepackage{amsfonts}
\usepackage{amssymb}
\usepackage{enumitem}
\usepackage{graphicx}
\usepackage{psfrag}
\usepackage{tikz}
\usepackage{subfigure}
\usepackage{makecell}
\usepackage{setspace}
\usepackage{scalerel}
\usepackage{booktabs}
\usepackage{authblk}

\counterwithin{figure}{section}
\counterwithin{table}{section}

\newlength{\blank}
\newtheorem{theorem}{Theorem}[section]
\newtheorem{lemma}[theorem]{Lemma}
\newtheorem{proposition}[theorem]{Proposition}
\newtheorem{corollary}[theorem]{Corollary}
\newtheorem{conjecture}[theorem]{Conjecture}

\newtheorem{statement}{\hspace{-\blank}}[theorem] 

\theoremstyle{definition}

\newcommand{\breaksta}{\hfill ~\linebreak }

\DeclareMathOperator{\cross}{cr}

\DeclareMathOperator{\sign}{sign}
\DeclareMathOperator{\boldsign}{\mathbf{sign}}

\DeclareMathOperator{\pf}{pf}
\DeclareMathOperator{\spf}{spf}

\DeclareMathOperator{\Pf}{Pf}
\DeclareMathOperator{\rank}{rank}

\DeclareMathOperator{\rev}{\mathbin{rev}}

\DeclareMathOperator{\genus}{g}
\DeclareMathOperator{\megenus}{\genus_\pmset}

\newcommand{\sharpp}{{\footnotesize\#}\textsc{p}}
\DeclareMathOperator{\bigO}{\mathcal{O}}

\newcommand{\dotprod}[2]{\langle #1, #2 \rangle}

\newcommand{\symdiff}{\mathbin{\triangle}}
\newcommand{\pmset}{\mathcal{M}}

\newcommand{\gf}{\textnormal{GF}}

\renewcommand{\vec}[1]{\mathbf {#1}}

\title{
The symbolic pfaffian method I: A generalization of the Gallucio-Loebl-Tesler Theorem
	\footnote{
		This study was financed, in part, by the São Paulo Research Foundation (FAPESP), Brasil. Process Number \#2026/00996-9. \\
		E-mail address: enriquejh@usp.br (E. Junchaya)
	}
}
\author[a]{Enrique Junchaya}
\author[a]{Cláudio L. Lucchesi}
\author[b]{Alberto A. A. Miranda}
\affil[a]{Universidade de São Paulo, Brazil 05508-020}
\affil[b]{Instituto Federal do Norte de Minas Gerais, Montes Claros, Brazil 39404-058}

\date{}

\begin{document}
	\maketitle

\begin{abstract}
	The pfaffian number of a graph $G$, denoted by $\pf(G)$, is the minimum integer~$k$ such that the generating function of perfect matchings of $G$ can be written as a linear combination of $k$ pfaffians of signed adjacency matrices.
	The Galluccio-Loebl-Tesler Theorem states that $\pf(G) \leq 4^{\genus(G)}$, where $\genus(G)$ is the orientable genus of $G$.
	Norine strengthened this result to $\pf(G) \leq 4^{\megenus(G)}$, where $\megenus(G)$ is the minimum genus of an orientable surface on which $G$ can be drawn such that every perfect matching of $G$ has an even number of crossings.

	As a consequence, $\pf(K_{2n}) = 2^{\bigO(n^2)}$ since the genus of complete graphs grows quadratically. 
	We improve this upper bound to $\pf(K_{2n}) = 2^{\bigO(n \log n)}$ by showing that $\pf(G) \leq |\pmset(G)|$, where $\pmset(G)$ is the set of perfect matchings of~$G$.

	Our main result is that $\pf(G) \leq 2^{\spf(G)}$, where $\spf(G)$ is the symbolic pfaffian number of $G$, another graph parameter.
	This upper bound unifies the previous known bounds.
	It is known that $\spf(G) = \bigO(\log |\pmset(G)|)$, from which one can obtain $\pf(G) \leq |\pmset(G)|^c$, with $c$ constant, which still implies $\pf(K_{2n}) = 2^{\bigO(n \log n)}$.
	We also prove that $\spf(G) \leq 2\megenus(G)$, which implies Norine's bound.

	Furthermore, we consider lower bounds for $\pf(G)$ related to $\spf(G)$.
	In general, we have that $\spf(G) + 1 \leq \pf(G)$, but we improve it in some cases;
	for instance, we show that $\pf(G) = 4$ if and only if $\spf(G) = 2$.
\end{abstract}

\tableofcontents

\section{Introduction}

This article is concerned with pfaffian-based methods for counting perfect matchings in finite graphs.
Hence, edges that do not belong to any perfect matching are not relevant.
Moreover, the number of perfect matchings of a disconnected graph is equal to the product of the number of the perfect matchings of its components.
Thus we can naturally restrict ourselves to matching covered graphs.
A graph $G$ is \emph{matching covered} if it is connected, has at least one edge and every edge of $G$ is contained in a perfect matching of $G$.
Throughout this work, whenever we say graph, assume we mean matching covered graph.

Assume that the graphs in this work are finite and simple.
Let~$G$ be a graph and let $x_e$, for every $e \in E(G)$, be a variable.
Let $\mathbf{x} = (x_e : e \in E(G))$ be a vector.
We define \emph{the generating function of perfect matchings} of $G$ by
\begin{equation}\label{eq:generating-function}
	\pmset_G(\mathbf{x}) := \sum_{M \in \pmset(G)} \prod_{e \in M} x_e,
\end{equation}
where $\pmset(G)$ is the set of perfect matchings of $G$.
Notice that, by evaluating the generating function of perfect matchings in $\mathbf{1}$ (the all-ones vector), we have
	\[ \pmset_G(\mathbf{1}) = |\pmset(G)|, \]
the number of perfect matchings of $G$.
In a more general setting, we can deal with graphs with weights on their edges.
The \emph{weight} of a perfect matching is equal to the product of the weights of its edges.
Let $\vec w := (w_e : e \in E(G))$ be the vector of weights of the edges of $G$.
In that case, $\pmset_G(\vec w)$ would be equal to the sum of the weights of all perfect matchings of~$G$.
Notice that multigraphs can be modeled as weighted simple graphs with weights accounting for multiplicity of the edges, thus we can assume that graphs are simple.
The (weighted) enumeration of perfect matchings is an interesting problem from a theoretical point of view, but it also has been studied in applied sciences.
The main applications are in chemistry (see \mbox{Gutman}~(1977, \cite{gutm1977}) and \mbox{Herndon}~(1974, \cite{hern1974})), statistical physics (see Brush~(1967, \cite{brus1967}) and \mbox{Kasteleyn}~(1967,~\cite{kast1967})) and quantum mechanics (see Krem, Gu and Zeilinger~(2017, \cite{kgz2017})).
Examples of some of these applications can be found in the book by Lovász and Plummer~(1986,~\cite[Chapter~8]{lopl1986}).
See also the lecture notes of Kenyon~(2009,~\cite{keny2009}).

For simplicity, we will restrict ourselves to the problem of counting perfect matchings, i.e. the case in which every edge has weight equal to one.
However, it is not hard to see that the arguments on this work are easily generalized to graphs with weighted edges. 

In bipartite graphs, counting perfect matchings is equivalent to computing the permanent of a square matrix containing only zeros and ones.
Valiant~(1979,~\cite{vali79}) proved that this problem is \mbox{\sharpp-complete};
therefore, the algorithms discussed here are not efficient in general.
In fact, the best known algorithm for counting perfect matchings is due to Li~(2026,~\cite{li2026}) and runs in $\bigO^*(2^{n/2 - \Omega(\sqrt n)})$ for a graph $G$ with $n = |V(G)|$ ($\bigO^*$ ignores factors that are polynomial on the input size).
It is worth mentioning that Li improves on a previous algorithm by Bj\"orklund~(2012,~\cite{bjor2012}).
In this paper, we focus on algorithms for counting perfect matchings using pfaffians.
Li's algorithm has better run time than these algorithms in general, but the pfaffian-based methods perform better in some particular graph classes, especially the graphs with small genus.

\subsection{Pfaffians and perfect matchings}

Let $G$ be a graph and let $D$ be an orientation of $G$.
Let $(1, 2, \dots, 2k)$ be an enumeration of the vertices in $G$.
Let $M$ be a perfect matching of $G$.
The \emph{sign} of $M$ in $D$, denoted by $\sign_D(M)$, is the sign of the permutation
	\[
		\begin{pmatrix}
			1 & 2 & 3 & 4 & \dots & 2k -1 & 2k \\
			u_1 & v_1 & u_2 & v_2 & \dots & u_k & v_k \\
		\end{pmatrix},
	\]
where, for each $i$, with $1 \leq i \leq k$, $u_iv_i$ is an edge of $M$ oriented from $u_i$ to $v_i$ in~$D$.
Notice that $\sign_D(M)$ is independent of the order in which the edges of $M$ are listed in the permutation.
Also, notice that reordering the vertices of $G$ either changes the sign of all of its perfect matchings or does not change the sign of any of its perfect matchings.
Assume that each graph has a fixed enumeration of its vertices, hence it makes sense to compare the signs of a perfect matching $M$ of $G$ in two distinct orientations of $G$.

Denote the adjacency matrix of the orientation~$D$ by $\mathbf{A}_D := (a_{ij})$, defined by
	\[
		a_{ij} := \left\{
		\begin{array}{rl}
			1,  & \text{if $ij$ is an arc of } D, \\
			-1, & \text{if $ji$ is an arc of } D, \\
			0,  & \text{otherwise}. \\
		\end{array}
		\right.
	\]
for $1 \leq i, j \leq |V(G)|$.
Notice that this matrix is skew-symmetric.
The \emph{pfaffian} of a skew-symmetric matrix~$\mathbf{A}$, denoted by $\Pf(\mathbf{A})$, is a polynomial in the entries of $\mathbf{A}$.
The pfaffian of a matrix can be computed efficiently; in fact, its square is equal to the determinant of the same matrix (Cayley, (1849,~\cite{cayl1849})).
For adjacency matrices of orientations, we can write the definition of the pfaffian as
	\[ \Pf(\mathbf{A}_D) := \sum_{M \in \pmset(G)} \sign_D(M). \]
An orientation $D$ of $G$ is \emph{pfaffian} if all perfect matchings of $G$ have the same sign in~$D$.
A graph is \emph{pfaffian} if it admits a pfaffian orientation.
Notice that, for a pfaffian orientation $D$, $|\Pf(\vec A_D)|$ is equal to the number of perfect matchings (and thus this number is efficiently computable).
Kasteleyn~(1963,~\cite{kast63}), who introduced pfaffian orientations, proved that every planar graph is pfaffian.

There is a vast theory and research on pfaffian graphs and orientations.	
Probably, the most astounding result is the efficient algorithm for deciding if a bipartite graph is pfaffian due to Robertson, Seymour and Thomas~(1999, \cite{rst1999}), and independently McCuaig~(2004, \cite{mccu2004}).
We refer the reader to the book of Lucchesi and Murty~(2024, \cite[Part III]{lumu2024}) for a detailed account of the main results on the theory of pfaffian graphs.

There are generalizations of the pfaffian method to count perfect matchings discussed above.
First, we introduce the $k$-pfaffian method.
A \emph{$k$-orientation} $\vec D$ of a graph $G$ is a vector $(D_1, D_2, \dots, D_k)$ of orientations of $G$.
The \emph{sign} of a perfect matching $M$ in $\vec D$ is defined as
	\[ \boldsign_{\vec D}(M) := (\sign_{D_1}(M), \sign_{D_2}(M), \dots, \sign_{D_k}(M)). \]
The \emph{signature matrix} $\boldsign_{\vec D}$ of $\vec D$ is a matrix with $|\pmset(G)$| rows and $k$ columns, in which every row is the sign in $\vec D$ of a different perfect matching of $G$.
The \mbox{$k$-orientation} $\vec D$ is \emph{pfaffian} if the system
	\[ \boldsign_{\vec D} \vec x = \vec 1 \]
has a solution.
In that case, we call that solution a \emph{solution} of $\vec D$.
A graph is \mbox{\emph{$k$-pfaffian}} if it admits a pfaffian $k$-orientation.
In particular, the $1$-pfaffian graphs are the pfaffian graphs.
Given a pfaffian $k$-orientation $\vec D$ of $G$ and a solution $\vec x$, one can count the number of perfect matchings of $G$ by the linear combination
	\[ \sum_{i = 1}^k x_i \Pf(\vec A_{D_i}). \]

The \emph{(orientable) genus} of a graph~$G$, denoted by $\genus(G)$, is the minimum integer $g$ such that $G$ is embeddable on the orientable surface of genus~$g$.
Kasteleyn~(1967, \cite{kast1967}) claimed, without proof, that it is possible to compute the number of perfect matchings of a graph~$G$ as a linear combination of~$4^{\genus(G)}$ pfaffians of oriented adjacency matrices.
This claim was formally proved by the end of the 20th century by Galluccio and Loebl~(1998,~\cite{galo98a}), and independently by Tesler~(2000,~\cite{tesl00}).
The \emph{pfaffian number} of a graph $G$, denoted by~$\pf(G)$, is the minimum integer $k$ such that $G$ is $k$-pfaffian.
Then, the Galluccio-Loebl-Tesler Theorem can be stated as follows.
\begin{theorem}[\cite{galo98a, tesl00}]\label{thm:glt}
	Let $G$ be a graph. Then,
		\[ \pf(G) \leq 4^{\genus(G)}. \eqno{\qed} \]
\end{theorem}

Norine~(2009,~\cite{nori09}) generalized the Galluccio-Loebl-Tesler Theorem by considering drawings of graphs with (restricted) intersections.
Let $\Sigma$ be a drawing of $G$ (possibly with crossings) on an orientable surface.
Denote by $\cross_\Sigma(M)$ the number of crossings of a perfect matching $M$ of $G$ in the drawing $\Sigma$.
The drawing $\Sigma$ is \emph{matching-even} if~$\cross_\Sigma(M)$ is even for every perfect matching $M$ of $G$.
The \emph{matching-even genus} of $G$, denoted by $\megenus(G)$, is the minimum integer $g$ such that $G$ admits a matching-even drawing on a surface of genus $g$.
\begin{theorem}[\cite{nori09}]\label{thm:norine-megenus}
	Let $G$ be a graph. Then,
		\[ \pf(G) \leq 4^{\megenus(G)}. \eqno{\qed} \]
\end{theorem}

In fact, equality holds in the theorem above for $\megenus(G) \leq 1$.
\begin{theorem}[\cite{nori09}]\label{thm:norine-megenus-equality}
	Let $G$ be a graph. Then,
		\[ \pf(G) = 1 \iff \megenus(G) = 0 \text{ and } \pf(G) = 4 \iff \megenus(G) = 1. \eqno{\qed} \]
\end{theorem}

Norine also showed that there are no graphs with pfaffian number equal to two, three or five.
\begin{theorem}[\cite{nori09}]\label{thm:non-existent-pf}
	Let $G$ be a graph. 
	If $\pf(G) \leq 3$, then $\pf(G) = 1$.
	If $\pf(G) \leq 5$, then $\pf(G) \leq 4$.
	\qed
\end{theorem}

It is also worth mentioning that there are known lower bounds for the pfaffian number.
In an earlier work, we proved that the pfaffian number is unbounded by showing that, for some infinite families of graphs, the pfaffian number grows at least linearly~(2026,~\cite{jml2026a}).
Then, Pant and Singh~(2026, \cite{pasi2026}) greatly improved this result by obtaining exponential lower bounds.
In particular, for complete graphs they obtained the following result.
\begin{theorem}\label{thm:pant-singh}
	For every $n \geq 3$,
		\[ \pf(K_{2n}) \geq \left\lceil \left( \frac{8}{3} \right)^{\lfloor n/3 \rfloor} \right\rceil. \]
\end{theorem}

Now we introduce yet another generalization of the pfaffian method.

\subsection{The symbolic pfaffian method}

Tesler~(2000,~\cite{tesl00}) pointed out that his proof of Theorem~\ref{thm:glt} also allowed counting perfect matchings by calculating the pfaffian of a single symbolic matrix.
This observation by Tesler gave rise to the symbolic pfaffian method.
Later, Babenko and Vyalyi~(2017,~\cite{bavy2017}) and Vyalyi~(2021,~\cite{vyal2021}) formalized this method.
The notations and definitions we use for this method are new, but based on the cited works.

We will denote $\gf(2)$ by $\mathbb{F}$ throughout this text.
Let $G$ be a graph and let $d$ be an integer.
A \emph{$d$-symbolic function}~$\tau$ of $G$ is a function whose domain is $E(G)$ and whose codomain is $\mathbb{F}^d$.
The function~$\tau$ extends to every subset of edges $S \subseteq E(G)$ by
	$$\tau(S) := \sum_{e \in S} \tau(e).$$
In particular, $\tau$ is applied to perfect matchings.
A pair $(D, \tau)$ is a \emph{$d$-symbolic orientation} of~$G$ if $D$ is an orientation of~$G$ and~$\tau$ is a $d$-symbolic function of $G$.
A $d$-symbolic orientation $(D, \tau)$ is $\emph{pfaffian}$ if for every pair of matchings $M, N \in \mathcal{M}(G)$ we have
	$$\sign_D(M) \neq \sign_D(N) \implies \tau(M) \neq \tau(N).$$
In other words, $(D, \tau)$ is pfaffian if $\tau$ separates perfect matchings of different signs.

Vyalyi~\cite{vyal2021} described an algorithm for counting the number of perfect matchings of $G$ given a pfaffian $d$-symbolic orientation of $G$ (see also Tesler~\cite[Appendix A.4.]{tesl00}).
Here, we outline this algorithm.
Let $p$ be the function that maps every vector $\mathbf{v}$ in~$\mathbb{F}^d$ to a monomial in the variables $t_i$, for $1 \leq i \leq d$, such that the exponent of $t_i$ in $p(\mathbf{v})$ is $v_i$.
Formally,
	\[
		\begin{array}{lccc}
			p: 	& \mathbb{F}^d & \longrightarrow{} & \mathcal{R}_d \\
				& \mathbf{v} & \longmapsto & \prod_{i = 1}^d t_i^{v_i}
		\end{array}
	\]
where $\mathcal{R}_d := \mathbb{Z}[t_1, \dots, t_d]/(t_i^2 - 1 : 1 \leq i \leq d)$ is the ring of polynomials in~$d$ variables with binary exponents.
Due to the restriction $t_i^2 = 1$, we have $p(\vec v + \vec w) = p(\vec v) \cdot p(\vec w)$ for all $\vec v, \vec w \in \mathbb{F}^d$.
Let $(D, \tau)$ be a pfaffian $d$-symbolic orientation of a graph $G$.
Let $\vec A_{(D,\tau)} = (a_{ij})$, for $1 \leq i, j \leq |V(G)|$, be a skew-symmetric matrix obtained by multiplying the entries of $\vec A_D$ by the corresponding monomial $p(\tau(e))$.
Formally,
	\[
		a_{ij} := \left\{
		\begin{array}{rl}
			p(\tau(ij)),  & \text{if $ij$ is an arc of } D, \\
			-p(\tau(ij)), & \text{if $ji$ is an arc of } D, \\
			0,  & \text{otherwise}. \\
		\end{array}
		\right.
	\]
By definition of $p$, we have $\prod_{e \in M} p(\tau(e)) = p(\tau(M))$, for every $M \in \pmset(G)$.
Therefore,
	\[	\Pf(\vec A_{(D, \tau)}) = \sum_{M \in \mathcal{M}(G)} \sign_D(M) p(\tau(M)), \]
a polynomial on the variables $t_i$.
By the definition of a pfaffian $d$-symbolic orientation, we know that the terms on the polynomial above do not cancel out.
In particular, the absolute value of the coefficient of a fixed monomial is equal to the number of perfect matchings $M$ that have $p(\tau(M))$ equal to that monomial.
Therefore, the number of perfect matchings of $G$ can be computed as the sum of the absolute values of the coefficients of the polynomial $\Pf(\vec A_{(D,\tau)})$.
There is a dynamic programming algorithm by Mahajan, Subramanya and Vinay~(1999, \cite{msv1999}) for computing the pfaffian of a matrix on a polynomial number of operations and without divisions; hence, this algorithm can be applied to arbitrary rings.
Since the intermediate operations needed for computing $\Pf(\vec A_{(D,\tau)})$ are performed over polynomials with at most $2^d$ terms, the complexity of this algorithm for counting perfect matchings is $\bigO^*(2^d)$.

Every graph $G$ admits a pfaffian $d$-symbolic orientation for sufficiently large $d$, as we will see.
The \emph{symbolic pfaffian number} of~$G$, denoted by $\spf(G)$, is the smallest integer $d$ such that the graph~$G$ admits a pfaffian $d$-symbolic orientation.
Notice that a pfaffian graph $G$ admits a pfaffian $0$-symbolic orientation $(D, \tau)$, where $D$ is a pfaffian orientation.
Thus, for a pfaffian graph $G$, $\spf(G) = 0$.

It may seem that a pfaffian $\spf(G)$-symbolic orientation $(D, \tau)$ of $G$ depends strongly on the orientation $D$.
However, for any orientation $D'$ of $G$ there is a $(\spf(G) + 1)$-symbolic function $\tau'$ such that $(D', \tau')$ is pfaffian.
Indeed, if $(D, \tau)$ is a pfaffian $\spf(G)$-symbolic orientation of $G$, then, for every edge $e \in E(G)$, it suffices to set $\tau'(e)$ equal to $\tau(e)$ in the first $\spf(G)$ entries and set the last entry of $\tau'(e)$ to be~$1$ if and only if the orientation of $e$ differs in $D'$ and $D$.
This last entry of $\tau'$ separates perfect matchings that have different signs in $D'$ but are mapped to the same vector by $\tau$.

\begin{proposition}
	Let $D$ be any orientation of a graph~$G$.
	There is a pfaffian $(\spf(G) + 1)$-symbolic orientation $(D, \tau)$.
	\qed
\end{proposition}

Vyalyi~\cite[proof of Lemma 1]{vyal2021} uses a standard probabilistic method to show that, for a set $S \subseteq \mathbb{F}^n$, there exists an injective linear map $f : S \rightarrow \mathbb{F}^d$ with $d = 2 \lceil \log_2 |S| \rceil$.
In other words, $f(x) \neq f(y)$, for distinct $x, y \in S$.
In particular, this result implies that there exists a $d$-symbolic function $\tau$ of~$G$, with $d$ equal to~$2 \lceil \log_2 |\pmset(G)| \rceil$, that maps different perfect matchings of $G$ to different vectors in $\mathbb{F}^d$.
\begin{theorem}\label{thm:spf-pmset-upper-bound}
	Let $G$ be a graph.
	Then,
		\[ \spf(G) \leq 2 \lceil \log_2 |\pmset(G)| \rceil. \eqno{\qed} \]
\end{theorem}

Since there are $2^{\bigO(n \log n)}$ perfect matchings in the complete graph, we obtain the following corollary.
\begin{corollary}\label{cor:spf-nlogn-upper-bound}
	Let $G$ be a graph with $n$ vertices. Then,
		\[ \spf(G) = \bigO(n \log n). \eqno{\qed} \]
\end{corollary}

This upper bound can be improved for graphs with small genus.
The following result is implicit in Tesler's proof of Theorem~\ref{thm:glt}.
\begin{theorem}\label{thm:spf-vs-genus}
	Let $G$ be a graph.
	Then,
		\[ \spf(G) \leq 2 \genus(G). \eqno{\qed} \]
\end{theorem}

Nevertheless, the upper bound in Corollary~\ref{cor:spf-nlogn-upper-bound} is tight.
\begin{theorem}[\cite{vyal2021}]
	$\spf(K_{2n}) = \Omega(n \log n)$ and $\spf(K_{n,n}) = \Omega(n \log n)$. \qed
\end{theorem}
Actually, Vyalyi proved the lower bounds above for a generalized version of the symbolic pfaffian method that allows multiple $d$-symbolic orientations.
These lower bounds hold for the symbolic pfaffian method as described here, a particular case of Vyalyi's generalization.

\subsection{Preliminaries and notation}

We use the graph theory notation of Bondy and Murty~(2008, \cite{bomu2008}).

\subsubsection*{Orientations and signs}

Let $G$ be a graph.
Let $D$ be an orientation of $G$ and  let $F$ be a set of arcs of $D$.
We denote by~$D \rev F$ the orientation obtained from $D$ by reversing the orientation of all the arcs of $F$.
The next proposition is a straightforward result.

\begin{proposition}\label{prop:rev-sign}
	Let $D$ be an orientation of $G$ and let $F$ be a set of arcs of $D$.
	Let~$M$ be a perfect matching of $G$.
	Then,
		\[ \sign_D(M) \cdot \sign_{D \rev F}(M) = (-1)^{|M \cap F|}. \eqno{\qed} \]
\end{proposition}

Let $D_1$ and $D_2$ be two orientations of the same graph.
We denote by $D_1 \symdiff D_2$ the set of edges whose orientations in $D_1$ and $D_2$ differ.

\subsubsection*{Linear algebra}

We denote by~$\vec J$ the square matrix with only $1$ entries and by $\vec I$ the identity matrix.
The dimensions of the square matrices $\vec J$ and $\vec I$ should be clear by context.

Let $\mathbf{u}$ and $\mathbf{v}$ be vectors in $\mathbb{R}^n$.
We denote the \emph{scalar product} of $\mathbf{u}$ and $\mathbf{v}$ by
	\[ \dotprod{\mathbf{u}}{\mathbf{v}} := \sum_{i = 1}^n u_i v_i. \]
We use this notation too for the scalar product of vectors in $\mathbb{F}^n$.
We say that $\mathbf{u}$ and~$\mathbf{v}$ are \emph{orthogonal} if $\dotprod{\mathbf{u}}{\mathbf{v}} = 0$.
Let $\mathcal{A}$ be a set of vectors in $\mathbb{R}^n$.
We say that $\mathcal{A}$ is \emph{orthogonal} if its elements are pairwise orthogonal.
The following proposition is a well-known property of orthogonal sets; see, for example, the book by Koffman and Hunze~(1971, \cite{hoku1971}).

\begin{proposition}\label{prop:orth->li}
	Let $\mathcal{A}$ be a set of non-zero vectors in $\mathbb{R}^n$.
	If $\mathcal{A}$ is orthogonal, then it is linearly independent.
	\qed
\end{proposition}

\subsection{Outline of the paper}

First, in Section~\ref{sec:pf-pmset-upperbound}, we prove that $\pf(G) \leq |\pmset(G)|$.
Given that the complete graph $K_{2n}$ has $2^{\bigO(n \log n)}$ perfect matchings, we obtain that 
\begin{equation}\label{eq:pf-nlogn}
	\pf(G) = 2^{\bigO(n \log n)},
\end{equation}
where $|V(G)| = 2n$, 
thus improving the previous known upper bound for complete graphs, i.e., $\pf(K_{2n}) = 2^{\bigO(n^2)}$, a consequence of the Galluccio-Loebl-Tesler Theorem and the fact that $\genus(K_{2n}) = \Theta(n^2)$ (Ringel and Youngs~(1968, \cite{riyo1968})).
Then, in Section~\ref{sec:main-theorem}, the main result is proved: the upper bound
	\[ \pf(G) \leq 2^{\spf(G)}. \]
This upper bound unifies the Galluccio-Loebl-Tesler Theorem and the bound in~\eqref{eq:pf-nlogn}.
Indeed, due to Theorem~\ref{thm:spf-vs-genus}, this result implies that $\pf(G) \leq 4^{\genus(G)}$.
Moreover, by Corollary~\ref{cor:spf-nlogn-upper-bound}, our main result implies that $\pf(K_{2n}) = 2^{\bigO(n \log n)}$.

In Section~\ref{sec:upperbounds-spf}, we derive some upper bounds for $\spf(G)$ based on $\pf(G)$.
In general, we have $\spf(G) \leq \pf(G) - 1$, but this can be improved in particular cases.
The most interesting result in that regard is that
	\[ \pf(G) = 4 \iff \spf(G) = 2. \]
In particular, by Theorem~\ref{thm:norine-megenus-equality}, we know that the graphs with $\spf(G) = 2$ are the graphs with $\megenus(G) = 1$.
More generally, in Section~\ref{sec:matching-even}, we improve Theorem~\ref{thm:spf-vs-genus} by considering matching-even drawings;
namely, we prove that
	\[ \spf(G) \leq 2 \megenus(G). \]
Notice that, if $\spf(G)$ were equal to $2 \megenus(G)$ in general, our main result would be equivalent to Theorem~\ref{thm:norine-megenus}.
We prove that this is not the case by showing a graph with odd symbolic pfaffian number in Section~\ref{sec:odd-spf}.

Finally, in Section~\ref{sec:conclusion}, we discuss conjectures and provide suggestions for further research.

\section{An algebraic upper bound for the pfaffian number}\label{sec:pf-pmset-upperbound}

There is an immediate algebraic upper bound for the pfaffian number that follows from the next known fact.

\begin{proposition}[\cite{jml2026a}]\label{prop:pf-rank}
	Let $\vec D$ be a pfaffian $k$-orientation of a graph $G$.
	Then,
		\[ \pf(G) \leq \rank(\boldsign_{\vec D}).  \eqno{\qed} \]
\end{proposition}

\begin{corollary}\label{cor:pf-pmset}
	For every graph $G$, we have
		\[ \pf(G) \leq |\pmset(G)|. \]
\end{corollary}
\begin{proof}
	Let $\vec D$ be a pfaffian $k$-orientation of $G$.
	By definition, $\boldsign_{\vec D}$ has at most $|\pmset(G)|$ rows.
	Thus, $\rank(\boldsign_{\vec D}) \leq |\pmset(G)|$.
	Then, by Proposition~\ref{prop:pf-rank}, we conclude that $\pf(G) \leq |\pmset(G)|$.
\end{proof}

Recall that $|\pmset(K_{2n})| = 2^{\bigO(n \log n)}$.
\begin{corollary}\label{cor:K2n-pf-upper-bound}
	Let $G$ be a graph with $2n$ vertices.
	Then, $\pf(G) = 2^{\bigO(n \log n)}.$ \qed
\end{corollary}

As the genus of the complete graph $K_{2n}$ is $\Theta(n^2)$, this improves the previous known upper bound of $\pf(K_{2n})$ based on the Galluccio-Loebl-Tesler Theorem, namely $\pf(K_{2n}) = 2^{\bigO(n^2)}$.
Likewise, $\genus(K_{n,n}) = \Theta(n^2)$~(Ringel~(1965,~\cite{ring1965})), thus Corollary~\ref{cor:K2n-pf-upper-bound} also improves the pfaffian number bound for complete bipartite graphs

On the other hand, notice that the upper bound of Corollary~\ref{cor:pf-pmset} only states the existence of a pfaffian $|\pmset(G)|$-orientation of a graph $G$, but does not imply an explicit construction of such orientation.

\section{A generalization of the Galluccio-Loebl-Tesler Theorem}\label{sec:main-theorem}

	We wish to prove that, for every graph $G$, $\pf(G) \leq 2^{\spf(G)}$ holds.
	In order to do that, we will construct a pfaffian $2^d$-orientation of a graph with $\spf(G) = d$.
	We will use the following lemma to show that the constructed $2^d$-orientation is pfaffian.

	\begin{lemma}\label{lem:pfaffian-half-difference-orientation}
		Let $\mathbf{D}$ be a $2k$-orientation of a graph.
		If every pair of different rows in~$\boldsign_\mathbf{D}$ differ in $k$ entries, then $\mathbf{D}$ is pfaffian.
	\end{lemma}
	\begin{proof}
		Let $\mathbf{A}$ be a matrix formed by the different rows of $\boldsign_\mathbf{D}$.
		If $\mathbf{A}$ has just one row, there is nothing to prove.
		Assume, then, that $\mathbf{A}$ has at least two rows.
		By hypothesis, every pair of different rows differ in half of the entries, then their scalar product is equal to zero.
		Thus, the rows of $\mathbf{A}$ form an orthogonal set and, by Proposition~\ref{prop:orth->li}, they are linearly independent.
		We conclude that the rank of~$\mathbf{A}$ is equal to the number of its rows, i.e., $\vec A$ is full rank.
		In particular, $\vec 1$ is in the column space of $\vec A$.
		In other words, the system $\vec A \vec x = \vec 1$ has solution.
		We conclude that $\vec D$ is pfaffian.
	\end{proof}

	To show that the condition of the lemma above holds for the $2^d$-orientation we will construct, we will use yet another lemma.
	\begin{lemma}\label{lem:half-of-binary-vectors}
		Let $\vec u \in \mathbb{F}^n$ be a non-zero vector.
		There are exactly $2^{n - 1}$ vectors $\vec k$ in~$\mathbb{F}^n$ such that
			\[ \dotprod{\vec u}{\vec k} \equiv 0 \pmod 2. \]
	\end{lemma}
	\begin{proof}
		By hypothesis, $\vec u \neq \vec 0$; thus, there is an integer $i$ such that $u_i = 1$.
		Let $\vec k \in \mathbb{F}^n$.
		Define $\overline{\vec k}$ as the vector in~$\mathbb{F}^n$ obtained from $\vec k$ by changing $k_i$.
		Clearly,
			\[ \dotprod{\vec u}{ \vec k} \not\equiv \dotprod{\vec u}{\overline{\vec k}} \pmod 2. \]
		Since $\vec k$ is arbitrary, half of the vectors $\vec k \in \mathbb{F}^n$ satisfy $\dotprod{\vec u}{\vec k} \equiv 0 \pmod 2$.
	\end{proof}

	\begin{theorem}\label{thm:main-theorem}
		Let $G$ be a graph.
		Then,
			$$ \pf(G) \leq 2^{\spf(G)}. $$
	\end{theorem}
	\begin{proof}
		Let $d := \spf(G)$.
		Let $(D, \tau)$ be a pfaffian $d$-symbolic orientation of $G$.
		We will show that $G$ is $2^d$-pfaffian by constructing a pfaffian $2^d$-orientation of this graph.
		For every vector $\vec i \in \mathbb{F}^d$, define the set of edges
		\[
			S_{\vec i} := \{e \in E(G) : \dotprod{\tau(e)}{\vec i} \equiv 1 \pmod 2\}.
		\]
		In particular, we have $S_{\vec 0} = \varnothing$.
		Let $X \subseteq E(G)$ be a set of edges.
		By the linearity of the scalar product and the definition of $S_{\vec i}$, we have
		\begin{equation}\label{eq:tau-X-dot-i}
			\dotprod{\tau(X)}{\vec i} \equiv \sum_{e \in X} \dotprod{\tau(e)}{\vec i} \equiv \sum_{e \in X \cap S_{\vec i}} 1 \equiv |X \cap S_{\vec i}| \pmod 2.
		\end{equation}

		Let $\mathbf{D} := (D_0, D_1, \dots, D_{2^d - 1})$ be a $2^d$-orientation of $G$ defined by
		\begin{equation}\label{eq:Di-definition}
			D_i := D \rev S_{\vec i}, 
		\end{equation}
		where $\vec i$ is the binary representation of the integer $i$ with $d$ bits.
		In particular, we have $D_0 = D$.
		By Proposition~\ref{prop:rev-sign}, and equations~\eqref{eq:tau-X-dot-i} and~\eqref{eq:Di-definition}, we conclude that
		\begin{equation}\label{eq:m-sign-in-Di}
			 \sign_{D_i}(M) \cdot \sign_{D}(M) = (-1)^{|M \cap S_{\vec i}|} = (-1)^{\dotprod{\tau(M)}{\vec i}}.
		\end{equation}
		We will show that $\mathbf{D}$ is pfaffian.

		If $d = 0$, then all perfect matchings of $G$ in $D$ have the same sign because $(D, \tau)$ is pfaffian.
		Therefore, $G$ is $1$-pfaffian and the theorem follows trivially.
		Suppose, then, that~$d \geq 1$.
		Let $M$ and $N$ be perfect matchings of $G$.
		We want to show that the vectors $\boldsign_{\vec D}(M)$ and $\boldsign_{\vec D}(N)$ are either equal or they differ in half of the entries.
		If $\tau(M) = \tau(N)$, then $\sign_{D}(M) = \sign_{D}(N)$ because $(D, \tau)$ is pfaffian. 
		Moreover, by~\eqref{eq:m-sign-in-Di}, we know that $\boldsign_\mathbf{D}(M) = \boldsign_\mathbf{D}(N)$.
		Suppose, then,  that $\tau(M) \neq \tau(N)$.
		Let $\vec a := \boldsign_{\vec D}(M)$ and $\vec b := \boldsign_{\vec D}(N)$ be $0$-indexed vectors.
		Since $d \geq 1$, these vectors have an even number of entries.
		The positions in which $\mathbf{a}$ and $\mathbf{b}$ differ correspond to the integers $i$ for which $a_i b_i = -1$.
		For every $i \in \{0, \dots, 2^d -1\}$, we have
			\[ a_i b_i = \sign_{D}(M) \sign_{D}(N) (-1)^{\dotprod{\tau(M) + \tau(N)}{\vec i}} \]
		by~\eqref{eq:m-sign-in-Di} and the linearity of the scalar product.
		Since $\tau(M) \neq \tau(N)$, we know that the vector $\tau(M) + \tau(N)$ is non-zero.
		Thus, by Lemma~\ref{lem:half-of-binary-vectors}, there are exactly $2^{d - 1}$ positions~$i$ such that
			\[ a_i b_i = \sign_{D}(M) \sign_{D}(N). \]
		In the other half of the positions, $a_ib_i$ has the opposite sign.
		Therefore, $\mathbf{a}$ and $\mathbf{b}$ differ in half of the entries, independently of the values of $\sign_{D}(M)$ and $\sign_{D}(N)$.
		Since~$M$ and~$N$ are arbitrary perfect matchings of $G$ and by Lemma~\ref{lem:pfaffian-half-difference-orientation}, we conclude that~$\mathbf{D}$ is pfaffian and, hence, $\pf(G) \leq 2^d$, as we wanted to prove.
	\end{proof}

	The theorem above together with Theorem~\ref{thm:spf-vs-genus} implies the Galluccio-Loebl-Tesler Theorem.
	Actually, as we will see in Section~\ref{sec:matching-even}, the theorem above implies Norine's Theorem~\ref{thm:norine-megenus}.
	On the other hand, Theorem~\ref{thm:spf-pmset-upper-bound}, that states $\spf(G) \leq 2 \lceil \log_2 |\pmset(G)| \rceil$, implies that	$\pf(G) \leq 4 |\pmset(G)|^2$.
	Therefore, this is another way of proving Corollary~\ref{cor:K2n-pf-upper-bound}.
	In sum, Theorem~\ref{thm:main-theorem} unifies the known upper bounds for the pfaffian number.

	Besides, recall that, given a pfaffian $d$-symbolic orientation of a graph~$G$, we can count the number of perfect matchings of~$G$ in $\bigO^*(2^d)$ time complexity.
	Also, given a pfaffian \mbox{$k$-orientation} of $G$ and its solution, we can count the number of perfect matchings of~$G$ in  $\bigO^*(k)$ time complexity.
	The Theorem~\ref{thm:main-theorem} implies then that (asymptotically) it is always at least as good to find a pfaffian $\pf(G)$-orientation than to find a pfaffian \mbox{$\spf(G)$-symbolic} orientation.
	This, however, does not render the symbolic pfaffian method irrelevant, as we do not know how to compute neither $\pf(G)$ nor $\spf(G)$ efficiently.
	Indeed, it is possible that obtaining a pfaffian $\spf(G)$-symbolic orientation is easier than obtaining a pfaffian $\pf(G)$-orientation.

\section{Upper bounds for the symbolic pfaffian number}\label{sec:upperbounds-spf}

We proved an upper bound on $\pf(G)$ that depends on $\spf(G)$.
The pfaffian number~$\pf(G)$ can also be lower bounded by $\spf(G)$, as we showed in an earlier paper~(2026,~\cite{jml2026a}).
Implicitly in our proof, we used the following definition of the symbolic sibling of a $k$-orientation.

Let $\vec D$ be a $k$-orientation.
The \emph{symbolic sibling} of $\vec D$ is the $(k - 1)$-symbolic orientation $(D_k, \tau)$, where $\tau$ is defined by
	\[
		\tau(e)_i := 
		\left\{
		\begin{matrix}
			1, & \text{if } e \in D_i \symdiff D_k \\
			0, & \text{otherwise}
		\end{matrix}
		\right.
	\]
for every edge $e \in E(G)$ and every $i \in \{1, \dots, k - 1\}$.

\begin{lemma}\label{lem:symbolic-sibling}
	Let $\vec D$ be a $k$-orientation of a graph $G$.
	Let $(D_k, \tau)$ be the symbolic sibling of $\vec D$.
	Then, for every $i \in \{1, \dots, k - 1\}$ and every $M \in \pmset(G)$, we have
		\[ \tau(M)_i = 1 \iff \sign_{D_i}(M) \neq \sign_{D_k}(M). \]
\end{lemma}
\begin{proof}
	Let $M$ be a perfect matching of $G$.
	By the definition of symbolic sibling, we have that $\tau(M)_i = 1$ if and only if $|M \cap (D_i \symdiff D_k)|$ is odd.
	By Proposition~\ref{prop:rev-sign}, we know that $|M \cap (D_i \symdiff D_k)|$ is odd if and only if $\sign_{D_i}(M) \cdot \sign_{D_k}(M) = -1$.
	Hence, the lemma follows.
\end{proof}

In~\cite{jml2026a}, we proved that if $\vec D$ is pfaffian, then its symbolic sibling is also pfaffian.

\begin{theorem}[\cite{jml2026a}]\label{thm:spf-basic-upper-bound}
	Let $G$ be a graph.
	Then,
		\[ \spf(G) + 1 \leq \pf(G). \eqno{\qed} \]
\end{theorem}

This lower bound is sharp because, for a pfaffian graph $G$, we have $\pf(G) = 1$ and $\spf(G) = 0$.
Nevertheless, this lower bound can be improved in particular cases.
In this section, we will consider the cases of graphs with $\pf(G) = 4$ and graphs with~$\pf(G) = 6$.

\subsection{Pfaffian preserving transformations}

Let $(D, \tau)$ be a $d$-symbolic orientation of a graph $G$.
Let $d'$ be an integer and let
	\[ f : \mathbb{F}^d \rightarrow{} \mathbb{F}^{d'} \]
be a linear transformation.
Let $\tau' := f \circ \tau$.
For the $d'$-symbolic orientation $(D, \tau')$, we have
	\[ \tau'(M) = \sum_{e \in M} f(\tau(e)) = f(\tau(M)), \]
for every $M \in \pmset(G)$, since $f$ is a linear transformation.
If $(D, \tau)$ and $(D, \tau')$ are pfaffian, we say that $f$ is a \emph{pfaffian preserving $d'$-transformation} of $(D, \tau)$.

We will use a pfaffian preserving transformations to improve the upper bound for~$\spf(G)$ in graphs with~$\pf(G) = 4$ and~$\pf(G) = 6$.
The graphs with these pfaffian numbers are relevant because we know fundamental properties of them.
\begin{theorem}[Norine, (2009~\cite{nori09})]\label{thm:norine-characterization}
	Let $G$ be a graph with $\pf(G) = 4$.
	Then, $G$ admits a pfaffian $4$-orientation~$\mathbf{D}$ with $\boldsign_\mathbf{D} = \mathbf{J} - 2\mathbf{I}$, up to repeated rows.
	In particular, $\boldsign_{\vec D}(M)$ has exactly one negative entry for every $M \in \pmset(G)$.
	\qed	
\end{theorem}

Actually, the theorem above is a rewrite of Norine's result using our notation.

\begin{theorem}[Costa Moço, Miranda and Nunes da Silva~(2021, \cite{mmn2021})]\label{thm:6-orientation}
	Let $G$ be a graph with $\pf(G) = 6$.
	There exists a pfaffian $6$-orientation $\vec D$ of $G$ for which every row of~$\boldsign_{\vec D}$ has exactly two negative entries.
	\qed
\end{theorem}

We will prove a somewhat generic result that applies to graphs with pfaffian number equal to four or six.

\begin{proposition}
	Let $k$ and $c$ be integers with $c \leq k$.
	If a graph~$G$ admits a pfaffian $k$-orientation $\vec D$ such that every row of $\boldsign_{\vec D}$ has exactly $c$ negative entries, then
		\[ \spf(G) \leq 2c. \]
\end{proposition}
\begin{proof}
	By Theorem~\ref{thm:spf-basic-upper-bound}, we know that $\spf(G) \leq k - 1$, thus we can assume that $c < (k - 1)/2$.
	Let the $(k - 1)$-symbolic orientation $(D_k, \tau)$ be the symbolic sibling of~$\vec D$.
	By Lemma~\ref{lem:symbolic-sibling}, for every $M \in \pmset(G)$, if $\sign_{D_k}(M) = 1$, then $\tau(M)$ has exactly~$c$ non-zero entries;
	otherwise, $\tau(M)$ has exactly $(k - 1) - (c - 1) = k - c$ non-zero entries.
	Let $f : \mathbb{F}^{k - 1} \rightarrow \mathbb{F}^{2c}$ be the projection of vectors in $\mathbb{F}^{k - 1}$ to its first~$2c$ entries.
	We want to show that $f$ is a pfaffian preserving $2c$-transformation of~$(D_k, \tau)$.
	Let $M$ be a perfect matching of $G$.
	If $\sign_{D_k}(M) = 1$, then $f(\tau(M))$ has at most~$c$ non-zero entries. 
	Otherwise, if $\sign_{D_k}(M) = -1$, then projecting $\tau(M)$ to its first~$2c$ entries, loses at most $k - 1 - 2c > 0$ non-zero entries;
	thus $f(\tau(M))$ has at least $(k - c) - (k - 1 - 2c) = c + 1$ non-zero entries.
	Therefore, $f \circ \tau$ separates perfect matchings of different signs in $D_k$, hence $(D_k, f \circ \tau)$ is pfaffian.
	We conclude that $\spf(G) \leq 2c$.
\end{proof}

By Theorems \ref{thm:norine-characterization} and \ref{thm:6-orientation}, we can apply the proposition above to graphs~$G$ with~$\pf(G) = 4$ or $\pf(G) = 6$.
\begin{corollary}\label{cor:pf-4-or-6-spf}
	Let $G$ be a graph.
	Then,
		\[ \pf(G) = 4 \implies \spf(G) \leq 2 \text{ and } \pf(G) = 6 \implies \spf(G) \leq 4. \eqno{\qed}\]
\end{corollary}

\subsection{$4$-Pfaffian graphs}

There is a stronger result for the symbolic pfaffian number of graphs with $\pf(G) = 4$.
As we will prove, the graphs with $\pf(G) = 4$ are the graphs with $\spf(G) = 2$.
The only missing ingredient is to show that if a graph has $\spf(G) \leq 1$, then it is pfaffian.

\begin{proposition}\label{prop:spf-1}
	Let $G$ be a graph with $\spf(G) \leq 1$.
	Then, $\spf(G) = 0$.
\end{proposition}
\begin{proof}
	Let $(D, \tau)$ be a pfaffian $1$-symbolic orientation of $G$.
	If all perfect matchings in $D$ have the same sign, then $G$ is clearly pfaffian.
	Thus, assume that there are perfect matchings of $G$ with different signs in $D$.
	By definition of pfaffian $d$-symbolic orientations, perfect matchings with different signs in $D$ have different images in $\tau$;
	then, $\tau(M) = (0)$ if and only if $\sign_D(M) = c$, for a constant $c$.
	Let $F$ be the set of edges~$e$ of $G$ with $\tau(e) = (1)$.
	Let $D' := D \rev F$ be an orientation of $G$.
	By Proposition~\ref{prop:rev-sign}, we have that
		\[ \sign_D(M) \neq \sign_{D'}(M) \iff |M \cap F| \equiv 1 \pmod 2 \iff \tau(M) = (1) \]
	for every $M \in \pmset(G)$.
	Since $\tau(M) = (1)$ if and only if $\sign_D(M) = -c$, all perfect matchings have sign $c$ in $D'$.
	We conclude that $G$ is pfaffian and thus $\spf(G) = 0$.
\end{proof}

\begin{corollary}\label{cor:spf-2=pf-4}
	Let $G$ be a graph.
	Then,
		\[ \spf(G) = 2 \iff \pf(G) = 4. \]
\end{corollary}
\begin{proof}
\breaksta\noindent
($\Rightarrow$) If $\spf(G) = 2$, then, by Theorem~\ref{thm:main-theorem}, we have that $\pf(G) \leq 4$.
We obtain that $\pf(G) = 4$ by Theorem~\ref{thm:non-existent-pf} since $G$ is not pfaffian.

\noindent
($\Leftarrow$) If $\pf(G) = 4$, then, by Corollary~\ref{cor:pf-4-or-6-spf}, we have that $\spf(G) \leq 2$.
Due to Proposition~\ref{prop:spf-1} and since $G$ is not pfaffian, we conclude that $\spf(G) = 2$.
\end{proof}

\section{Matching-even drawings}\label{sec:matching-even}

In this section, we show that, for a graph $G$, $\spf(G) \leq 2 \megenus(G)$ holds, thus extending Theorem~\ref{thm:spf-vs-genus}.
Our proof is based on the proof of Theorem~\ref{thm:glt} by Tesler~(2000, \cite{tesl00}).
We rely on crossing orientations, whose existence was stated by Kasteleyn~(1967,~\cite{kast1967}).
Let~$G$ be a graph and let $\Sigma$ be a drawing of $G$ in the plane.
An orientation $D$ of $G$ is a \emph{crossing orientation} in $\Sigma$ if 
	\[ \sign_D(M) = (-1)^{\cross_\Sigma(M)} \]
for every perfect matching $M$ of $G$.
Tesler~\cite{tesl00} gave an explicit proof of the existence of crossing orientations.

\begin{theorem}\label{thm:crossing-orientation}
	Let $G$ be a graph and let~$\Sigma$ be a drawing of $G$ in the plane.
	There exists a crossing orientation of $G$ in~$\Sigma$.
	\qed
\end{theorem}

For finite graphs, their drawings can be restricted to \emph{compact surfaces}: bounded surfaces without $1$-dimensional boundaries.
So assume that every surface in this work is compact.

\begin{figure}[ht]
	\centering
	\includegraphics{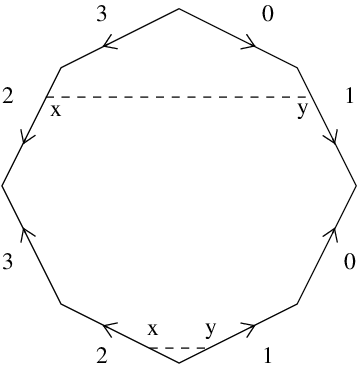}
	\caption{
		A plane model of the double torus, the orientable surface with genus equal to two.
		The points labeled $x$ represent a single point in the original surface.
		The same occurs with the points labeled $y$.
		The dotted line represents a curve on the surface.
		\label{fig:plane-double-torus}
	}
\end{figure}

It is convenient to represent non-planar surfaces in the plane.
We will work with a polygonal representation of surfaces, called the plane model of surfaces (see Firby and Gardiner~(2001,~\cite{figa2001})).
A \emph{plane model} of the orientable surface $S$ of genus $g$ is a solid polygon $P$ of $4g$ sides, drawn on the plane, with pairs of points in the sides identified as described below.
(See Figure~\ref{fig:plane-double-torus}, a plane model of the double torus.)
Enumerate the sides of $P$, in clockwise order, from $0$ to $4g - 1$.
For every $i \in \{0, \dots, g - 1\}$, the sides $4i$ and $4i + 2$ are labeled $a_{2i}$; also, the sides $4i + 1$ and $4i + 3$  are labeled~$a_{2i + 1}$.
Add arrows to the sides such that sides with the same label have arrows in opposite directions (considering the clockwise traversal of~$P$).
Consider the sides of $P$ that share the same label as identified, respecting the direction of the arrows.
By ``gluing'' the identified points of $P$ in a three-dimensional space, we obtain the original surface~$S$; thus, $P$ represents~$S$ faithfully.

\psfrag{0}{$a_0$}
\psfrag{1}{$a_1$}
\psfrag{2}{$a_2$}
\psfrag{3}{$a_3$}
\psfrag{x}{$x$}
\psfrag{y}{$y$}

\begin{figure}[ht]
	\centering
	\includegraphics[scale=1.25]{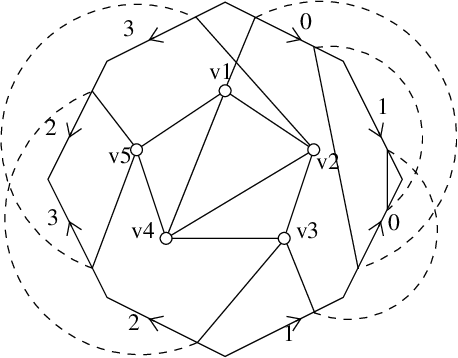}
	\caption{
		A planarization of a drawing of $K_5$ on the double torus.
		We did not try to minimize the number of crossings (in fact, $K_5$ is embeddable in the torus).
		\label{fig:planarization}
	}
\end{figure}

Now consider a drawing $\Sigma$ of a graph $G$ on the orientable surface of genus $g$.
A \emph{planarization} of $\Sigma$ is a drawing $\Sigma^P$ of $G$ in the plane obtained by taking $\Sigma$ into a plane model $P$ of $S$ and completing the edges that touch the sides of $P$ (on identified points) with curves (without self-intersections) outside of the polygon.
See Figure~\ref{fig:planarization}.
We say that a curve that connects two (identified) points of sides with the same label~$a_i$ is an \emph{$i$-curve}.
All \emph{$i$-curves}, for $i$ fixed, must be drawn in such a way that they do not cross each other.
They can be drawn as concentric arcs, for example.
Notice that every \mbox{$2i$-curve} crosses all the $(2i + 1)$-curves, for $i \in \{0, \dots, g - 1\}$.  
These curves must be drawn in a way such that each $2i$-curve crosses each \mbox{$(2i + 1)$-curve} once, and these crossings are the only crossings outside~$P$.
An edge may contain multiple $i$-curves.
Let $F$ be a set of edges of $G$.
Denote by~$N(F, i)$ the number of $i$-curves that are contained in the edges of~$F$ in the drawing~$\Sigma^P$.
For example, for the planarization in Figure~\ref{fig:planarization}, we have~$N(\{v_1 v_3\}, 0) = 2$ and~$N(\{v_1v_3\}, 1) = 1$.

\psfrag{v1}{$v_1$}
\psfrag{v2}{$v_2$}
\psfrag{v3}{$v_3$}
\psfrag{v4}{$v_4$}
\psfrag{v5}{$v_5$}

We will use the planarization of a matching-even drawing to prove the main result of this section.

\begin{theorem}\label{thm:spf-vs-megenus}
	Let $G$ be a graph.
	Then,
		\[ \spf(G) \leq 2 \megenus(G). \]
\end{theorem}
\begin{proof}
	Let $\megenus := \megenus(G)$.
	Let $S$ be the orientable surface of genus $\megenus$.
	Let $P$ be a plane model of $S$ and let $\Sigma$ be a planarization in $P$ of a matching-even drawing of $G$ on $S$.
	By Theorem~\ref{thm:crossing-orientation}, we know there is a crossing orientation $D$ of $G$ in the drawing~$\Sigma$.
	Define the $2\megenus$-symbolic function $\tau$ by
		\[
			\tau(e)_i := N(\{e\}, i) \pmod 2,
		\]	
	for $i \in \{0, \dots, 2\megenus - 1\}$.
	We want to show that $(D, \tau)$ is a pfaffian $2\megenus$-symbolic orientation.
	\begin{statement}
		For every $M \in \pmset(G)$, $\cross_{\Sigma}(M) \equiv \sum_{i = 0}^{\megenus - 1} \tau(M)_{2i} \cdot \tau(M)_{2i + 1} \pmod 2$.
	\end{statement}
	\begin{proof}
		Let $M$ be a perfect matching of $G$.
		Since $\Sigma$ is a planarization of a matching-even drawing, the number of crossings of $M$ in $\Sigma$ inside $P$ is even.
		Therefore only the crossings outside $P$ are relevant for counting $\cross_{\Sigma}(M)$ modulo~2.
		The only crossings outside $P$ are formed by $2i$-curves with $(2i + 1)$-curves, for every $i \in \{0, \dots, \megenus - 1\}$.
		Since each $2i$-curve crosses each $(2i + 1)$-curve once, the number of such crossings in~$M$ is $N(M, 2i) \cdot N(M, 2i + 1)$.
		Therefore,
		\begin{equation}\label{eq:sum-of-crossings}
			\cross_{\Sigma}(M) \equiv \sum_{i = 0}^{\megenus - 1} N(M, 2i) \cdot N(M, 2i + 1) \pmod 2.
		\end{equation}
		By the construction of $\tau$, we have that, for $j \in \{0, \dots, 2\megenus - 1\}$,
			\[N(M, j) \equiv \sum_{e \in M} N(\{e\}, j) \equiv \sum_{e \in M} \tau(e)_j \equiv \tau(M)_j \pmod 2 . \]
		The claim follows by replacing $N(M, j) \equiv \tau(M)_j \pmod 2$ in~\eqref{eq:sum-of-crossings}.
	\end{proof}
		
	Recall that, since $D$ is a crossing orientation in $\Sigma$, $\sign_D(M) = (-1)^{\cross_{\Sigma}(M)}$ for every perfect matching $M$ of $G$.
	Therefore, given that $\cross_{\Sigma}(M)$ is a function of~$\tau(M)$, we obtain that $\sign_D(M)$ is a function of $\tau(M)$.
	We conclude that, for every two perfect matchings $M_1$ and $M_2$ with different signs in~$D$, $\tau(M_1) \neq \tau(M_2)$ holds.
	By definition, $(D, \tau)$ is pfaffian and $\spf(G) \leq 2\megenus$, as we wanted to prove.
\end{proof}

Due to Theorem~\ref{thm:norine-megenus-equality} we know that the theorem above is sharp.
Indeed, the class of graphs with $\spf(G) = 2$ ($\pf(G) = 4$) is equal to the class of graphs with $\megenus(G) = 1$;
also, the graphs with $\spf(G) = 0$ are the graphs with $\megenus(G) = 0$.
We could not help but wonder if in general $\spf(G)$ is equal to $2\megenus(G)$.
This would imply that our main result, Theorem~\ref{thm:main-theorem}, is equivalent to Norine's Theorem~\ref{thm:norine-megenus}.
In the next section, we prove that this is not true by showing a graph with odd symbolic pfaffian number.

\section{A graph with odd symbolic pfaffian number}\label{sec:odd-spf}

\begin{figure}[ht]
	\centering
	\includegraphics[width=0.9\textwidth]{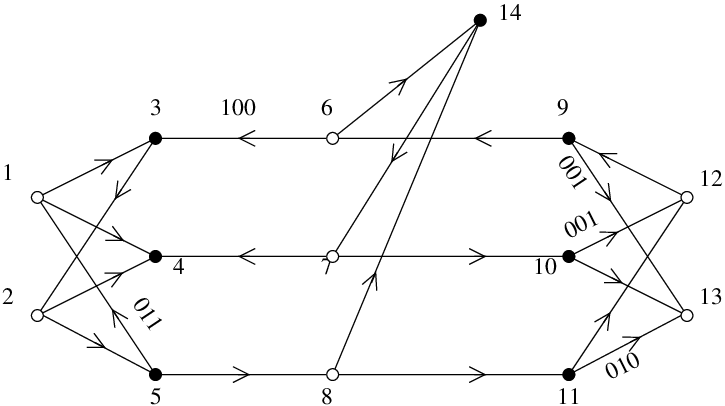}
	\caption{
		A pfaffian $3$-symbolic orientation $(D, \tau)$ of $G_{14,21,b}$.
		The vector $\tau(e)$ is shown at the side and parallel to the corresponding edge $e$.
		The zero vectors are omitted.
		\label{fig:spf-3}
	}
\end{figure}

In Figure~\ref{fig:spf-3} we have a pfaffian $3$-symbolic orientation of the graph $G_{14,21,b}$, a graph with symbolic pfaffian number equal to three.

\psfrag{1}{$v_1$}
\psfrag{2}{$v_2$}
\psfrag{3}{$v_3$}
\psfrag{4}{$v_4$}
\psfrag{5}{$v_5$}
\psfrag{6}{$v_6$}
\psfrag{7}{$v_7$}
\psfrag{8}{$v_8$}
\psfrag{9}{$v_9$}
\psfrag{10}{$v_{10}$}
\psfrag{11}{$v_{11}$}
\psfrag{12}{$v_{12}$}
\psfrag{13}{$v_{13}$}
\psfrag{14}{$v_{14}$}
\psfrag{100}{\scriptsize{$(1,0,0)$}}
\psfrag{011}{\scriptsize{$(0,1,1)$}}
\psfrag{010}{\scriptsize{$(0,1,0)$}}
\psfrag{001}{\scriptsize{$(0,0,1)$}}

To show that the $3$-symbolic orientation $(D, \tau)$ in Figure~\ref{fig:spf-3} is indeed pfaffian, we list all perfect matchings of $G_{14,21,b}$, their signs in $D$ and the vectors they are mapped to by $\tau$ in Table~\ref{tab:spf-3}.
It is clear that perfect matchings with different signs are separated by $\tau$; thus $(D, \tau)$ is pfaffian and $\spf(G_{14,21,b}) \leq 3$.
We also have a constructive proof to show that $\spf(G_{14,21,b}) \leq 3$, but we omit it here because the results necessary for that construction are out of the scope of this paper.
In an unpublished work, we proved that $\pf(G_{14,21,b}) > 4$ \cite{jml2026b}; hence, we can conclude that $\spf(G_{14,21,b}) = 3$ by Corollary~\ref{cor:spf-2=pf-4}.

\begin{table}[h]
	\centering
	\begin{tabular}{| c | c | c |}
		\toprule
		$M$ & $\sign_D(M)$ & $\tau(M)$ \\
		\midrule
		$\{v_{1}v_{3}, v_{2}v_{4}, v_{5}v_{8}, v_{6}v_{14}, v_{7}v_{10}, v_{9}v_{13}, v_{11}v_{12}\}$ & $+1$ & $(0,0,1)$ \\
		$\{v_{1}v_{3}, v_{2}v_{4}, v_{5}v_{8}, v_{6}v_{14}, v_{7}v_{10}, v_{12}v_{9}, v_{11}v_{13}\}$ & $+1$ & $(0,1,0)$ \\
		$\{v_{1}v_{3}, v_{2}v_{4}, v_{5}v_{8}, v_{9}v_{6}, v_{14}v_{7}, v_{10}v_{12}, v_{11}v_{13}\}$ & $-1$ & $(0,1,1)$ \\
		$\{v_{1}v_{3}, v_{2}v_{4}, v_{5}v_{8}, v_{9}v_{6}, v_{14}v_{7}, v_{10}v_{13}, v_{11}v_{12}\}$ & $+1$ & $(0,0,0)$ \\
		$\{v_{1}v_{3}, v_{2}v_{5}, v_{7}v_{4}, v_{6}v_{14}, v_{8}v_{11}, v_{9}v_{13}, v_{10}v_{12}\}$ & $+1$ & $(0,0,0)$ \\
		$\{v_{1}v_{3}, v_{2}v_{5}, v_{7}v_{4}, v_{6}v_{14}, v_{8}v_{11}, v_{12}v_{9}, v_{10}v_{13}\}$ & $+1$ & $(0,0,0)$ \\
		$\{v_{1}v_{3}, v_{2}v_{5}, v_{7}v_{4}, v_{9}v_{6}, v_{8}v_{14}, v_{10}v_{12}, v_{11}v_{13}\}$ & $-1$ & $(0,1,1)$ \\
		$\{v_{1}v_{3}, v_{2}v_{5}, v_{7}v_{4}, v_{9}v_{6}, v_{8}v_{14}, v_{10}v_{13}, v_{11}v_{12}\}$ & $+1$ & $(0,0,0)$ \\
		$\{v_{1}v_{4}, v_{2}v_{5}, v_{6}v_{3}, v_{7}v_{10}, v_{8}v_{14}, v_{9}v_{13}, v_{11}v_{12}\}$ & $+1$ & $(1,0,1)$ \\
		$\{v_{1}v_{4}, v_{2}v_{5}, v_{6}v_{3}, v_{7}v_{10}, v_{8}v_{14}, v_{12}v_{9}, v_{11}v_{13}\}$ & $+1$ & $(1,1,0)$ \\
		$\{v_{1}v_{4}, v_{2}v_{5}, v_{6}v_{3}, v_{14}v_{7}, v_{8}v_{11}, v_{9}v_{13}, v_{10}v_{12}\}$ & $-1$ & $(1,0,0)$ \\
		$\{v_{1}v_{4}, v_{2}v_{5}, v_{6}v_{3}, v_{14}v_{7}, v_{8}v_{11}, v_{12}v_{9}, v_{10}v_{13}\}$ & $-1$ & $(1,0,0)$ \\
		$\{v_{1}v_{4}, v_{3}v_{2}, v_{5}v_{8}, v_{6}v_{14}, v_{7}v_{10}, v_{9}v_{13}, v_{11}v_{12}\}$ & $+1$ & $(0,0,1)$ \\
		$\{v_{1}v_{4}, v_{3}v_{2}, v_{5}v_{8}, v_{6}v_{14}, v_{7}v_{10}, v_{12}v_{9}, v_{11}v_{13}\}$ & $+1$ & $(0,1,0)$ \\
		$\{v_{1}v_{4}, v_{3}v_{2}, v_{5}v_{8}, v_{9}v_{6}, v_{14}v_{7}, v_{10}v_{12}, v_{11}v_{13}\}$ & $-1$ & $(0,1,1)$ \\
		$\{v_{1}v_{4}, v_{3}v_{2}, v_{5}v_{8}, v_{9}v_{6}, v_{14}v_{7}, v_{10}v_{13}, v_{11}v_{12}\}$ & $+1$ & $(0,0,0)$ \\
		$\{v_{5}v_{1}, v_{2}v_{4}, v_{6}v_{3}, v_{7}v_{10}, v_{8}v_{14}, v_{9}v_{13}, v_{11}v_{12}\}$ & $+1$ & $(1,1,0)$ \\
		$\{v_{5}v_{1}, v_{2}v_{4}, v_{6}v_{3}, v_{7}v_{10}, v_{8}v_{14}, v_{12}v_{9}, v_{11}v_{13}\}$ & $+1$ & $(1,0,1)$ \\
		$\{v_{5}v_{1}, v_{2}v_{4}, v_{6}v_{3}, v_{14}v_{7}, v_{8}v_{11}, v_{9}v_{13}, v_{10}v_{12}\}$ & $-1$ & $(1,1,1)$ \\
		$\{v_{5}v_{1}, v_{2}v_{4}, v_{6}v_{3}, v_{14}v_{7}, v_{8}v_{11}, v_{12}v_{9}, v_{10}v_{13}\}$ & $-1$ & $(1,1,1)$ \\
		$\{v_{5}v_{1}, v_{3}v_{2}, v_{7}v_{4}, v_{6}v_{14}, v_{8}v_{11}, v_{9}v_{13}, v_{10}v_{12}\}$ & $-1$ & $(0,1,1)$ \\
		$\{v_{5}v_{1}, v_{3}v_{2}, v_{7}v_{4}, v_{6}v_{14}, v_{8}v_{11}, v_{12}v_{9}, v_{10}v_{13}\}$ & $-1$ & $(0,1,1)$ \\
		$\{v_{5}v_{1}, v_{3}v_{2}, v_{7}v_{4}, v_{9}v_{6}, v_{8}v_{14}, v_{10}v_{12}, v_{11}v_{13}\}$ & $+1$ & $(0,0,0)$ \\
		$\{v_{5}v_{1}, v_{3}v_{2}, v_{7}v_{4}, v_{9}v_{6}, v_{8}v_{14}, v_{10}v_{13}, v_{11}v_{12}\}$ & $-1$ & $(0,1,1)$ \\
		\bottomrule
	\end{tabular}
	\caption{
		Perfect matchings of $G_{14,21,b}$ with signs in $D$ and vectors in $\tau$ for the pfaffian $3$-symbolic orientation $(D, \tau)$ of Figure~\ref{fig:spf-3}
		\label{tab:spf-3}
	}
\end{table}

\section{Conclusion}\label{sec:conclusion}

	We have proven the following relations of the pfaffian number and the symbolic pfaffian number of a graph $G$:
		\[ \spf(G) + 1 \leq \pf(G) \leq 2^{\spf(G)}. \]
	The lower bound is sharp for pfaffian graphs, but we have seen that it can be improved for graphs with $\pf(G) = 4$ and $\pf(G) = 6$.
	We believe that the lower bound can be greatly improved in general.
	The fact that the cost for applying the symbolic pfaffian method to count the number of perfect matchings of $G$ is exponential in $\spf(G)$ suggests that the real lower bound for $\pf(G)$ is exponential in $\spf(G)$.	

	\begin{conjecture}
		$ \pf(G) = 2^{\Omega(\spf(G))}. $
	\end{conjecture}

	In particular, we know that $\spf(K_{2n}) = \Theta(n \log n)$ but $\pf(K_{2n})$ grows exponentially with $n$.
	Indeed, by Theorem~\ref{thm:pant-singh} and Corollary~\ref{cor:K2n-pf-upper-bound}, we know that
		\[ \left\lceil \left( \frac{8}{3} \right)^{\lfloor n/3 \rfloor} \right\rceil \leq \pf(K_{2n}) = 2^{\bigO(n \log n)}. \]
	It is an open problem to close this gap and settling our conjecture above would show that the upper bound is tight.
	
	On the other hand, we have concluded that the Galluccio-Loebl-Tesler Theorem, namely, the upper bound $\pf(G) \leq 4^{\genus(G)}$, is a particular case of a more general method that does not involve embeddings on surfaces.
	We do not know, however, how to efficiently compute $\spf(G)$, so our result does not directly imply an algorithm for counting perfect matchings.
	Despite that, since the proof of Theorem~\ref{thm:main-theorem} is constructive, it opens the possibility to obtain FPT algorithms by constructing pfaffian $d$-symbolic orientations with $d$ bounded by a graph parameter.

	It is worth mentioning that there are different versions of the symbolic pfaffian method, with vectors over abelian groups different from $\gf(2)$.
	For instance, Tesler~\cite{tesl00} considered the symbolic pfaffian method over $\gf(4)$ to compute the generating function of perfect matchings using embeddings on non-orientable surfaces.
	Also, Babenko and Vyalyi~\cite{bavy2017} and Vyalyi~\cite{vyal2021} found bounds for the symbolic pfaffian number over~$\mathbb{Z}$ (besides over $\gf(2)$).
	To the best of our knowledge, these other variants of the symbolic pfaffian method are largely unexplored.

	\section*{Acknowledgments} 
	We would like to thank Lucas Collucci for his comments that helped simplifying some of our linear algebra arguments.
	We are also grateful to Marcel K. de Carli Silva and Arnaldo Mandel for their comments on the writing of some excerpts of this text.

	\bibliographystyle{babplain}
	\bibliography{bibliografia}
\end{document}